\documentclass[english, a4paper, 10pt]{amsart}
\usepackage[utf8]{inputenc}
\usepackage{graphicx}
\usepackage{pgf,tikz}
\usetikzlibrary{arrows}

\usepackage[margin=2.7cm]{geometry}

\usepackage{amsthm,amsfonts,amssymb,amsmath,amscd}
\usepackage[all]{xy}
\usepackage{xr-hyper}
\usepackage{verbatim}
\usepackage{pdfsync}

\numberwithin{equation}{section}
\numberwithin{figure}{section}

\usepackage{babel}
\usepackage{ngerman}
\usepackage{graphicx}
\usepackage[T1]{fontenc}
\usepackage[utf8]{inputenc}
\usepackage{amsmath}
\usepackage{amssymb}
\usepackage{amsthm}
\usepackage{amsfonts}
\usepackage{amssymb}
\usepackage{hyperref}

\renewcommand{\omega}{\gamma}
\renewcommand{\Omega}{\Gamma}
\renewcommand{\nu}{\lambda}
\renewcommand{\rho}{\varrho}

\newcommand{\Adm}{\operatorname{Adm}}
\newcommand{\Perm}{\operatorname{Perm}}
\newcommand{\Wtilde}{\stackrel{\sim}{\smash{W}\rule{0pt}{1.1ex}}}
\newcommand{\inWtilde}{\in\ \Wtilde}
\newcommand{\Wmu}{W \! \mu}

\title{On the $\mu$-admissible set in the extended affine Weyl groups of $E_6$ and $E_7$}
\author{Lisa Sauermann}
\date{1st April 2016}

\begin{document}

\begin{abstract}Giving explicit examples we show that the $\mu$-admissible set and the $\mu$-permissible set do not agree for non-trivial minuscule coweights in the cases of $E_6$ and $E_7$.\end{abstract}

\maketitle

\section{Introduction}
In \cite{kotrapo},  Kottwitz and Rapoport introduce the notions of the $\mu$-admissible set $\Adm(\mu)$ and $\mu$-permissible set $\Perm(\mu)$ for a cocharacter $\mu$ (we will recall these notions in section \ref{not}) and in \cite[\S11]{kotrapo} they prove that $\Adm (\mu)\subset \Perm (\mu)$. Haines and Ngo show in \cite[Theorem 3]{hn}, that the equality $\Adm (\mu)=\Perm (\mu)$ does not hold for general $\mu\in X_{*}$.  On the other hand, Kottwitz and Rapoport prove in \cite{kotrapo} that $\Adm (\mu)= \Perm (\mu)$ holds for \emph{minuscule} cocharacters $\mu$ for the root systems $A_n$ and $C_n$.  They raise the question whether the $\mu$-permissible set $\Perm(\mu)$ and the $\mu$-admissible set $\Adm(\mu)$ agree for any minuscule cocharacter $\mu$, comp. also \cite[after Thm. 3.4]{rapoguide}. Smithling proves  in \cite{sm} that  $\Adm (\mu)=\Perm (\mu)$ for minuscule cocharacters $\mu$ for the root systems $B_n$ and $D_n$. Thus, this question has a positive answer  for all classical root systems.

The aim of the present paper is to show that the question has a negative answer in the cases  $E_6$ and $E_7$. We follow the following strategy. 

In section \ref{e6} we consider a certain root datum for the root system $E_6$ and a certain minuscule cocharacter $\mu$. We exhibit  a certain element $x$ of its extended affine Weyl group (everything will be defined in subsection \ref{settin6}), and investigate whether $x$ is $\mu$-permissible, resp. $\mu$-admissible. In subsection \ref{xperm6} we show that $x$ is $\mu$-permissible using  only simple matrix calculations. Afterwards, in subsection \ref{xadm16}, we find that $x$ is not $\mu$-admissible relying on a result of He and Lam which characterizes the $\mu$-admissible set. As a double-check, we prove in subsection \ref{xadm26} that $x$ is not $\mu$-admissible using computer software, this time relying on a result of Haines, instead of the characterization by He and Lam. Thus, the existence of $x$ gives a negative answer to the question of Kottwitz and Rapoport. 

In section \ref{e7} we repeat the same considerations for a certain element of the extended affine Weyl group attached to some root datum for the root system $E_7$ and some minuscule coweight $\mu$, again yielding a counterexample to the question of Kottwitz and Rapoport.

The counterexamples presented in this paper were found using the CHEVIE-package of the software GAP. The author wrote a program that calculated the size of the permissible and the admissible set in the above settings. This gave 20159 for the cardinality of the $\mu$-admissible set and 20303 for the cardinality of the $\mu$-permissible set for $E_6$; for $E_7$, the $\mu$-admissible set  has 1227151 elements and the $\mu$-permissible set  has 1298607 elements.

As the conjecture of Kottwitz and Rapoport does not depend on the actual choice of the root datum, our calculations show that $\Adm(\mu)\neq \Perm(\mu)$ holds in general for the considered minuscule coweights of $E_6$ and $E_7$. Note that $E_7$ has only one non-trivial dominant minuscule coweight, and in $E_6$ the two non-trivial dominant minuscule coweights are interchanged by an automorphism of the root system of $E_6$. Hence we get $\Adm(\mu)\neq \Perm(\mu)$ for all  non-trivial minuscule coweights $\mu$ of $E_6$ and $E_7$. As $\Adm(\mu)=\Perm(\mu)$ has been proved for minuscule coweights $\mu$  for classical root systems, this finally gives a complete  answer to the question on the relation between the $\mu$-admissible and the $\mu$-permissible set for minuscule coweights $\mu$ (note that in the other cases $E_8$, $F_4$ and $G_2$ there are no non-trivial minuscule coweights).

In the first version of this paper the author checked $x\not\in \Adm(\mu)$ only by computer calculations. The referee suggested how to prove $w_2\not\leq w_1$ in subsections \ref{xadm16} and \ref{xadm17} by hand, eliminating the need of computer calculations. However, in subsection \ref{xadm26} we still present an alternative way for checking $x\not\in \Adm(\mu)$ using the computer software Sage, that does not rely on the characterization of the $\mu$-admissible set by He and Lam.

The author wants to thank Michael Rapoport for his continuous support. Without his encouragement and advice this work would have never been possible. Furthermore the author wants to thank Xuhua He for a very helpful conversation, where he pointed out his criterion used in subsections \ref{xadm16} and \ref{xadm17}, as well as Ulrich Görtz and Robert Kottwitz for their helpful comments on earlier versions of this paper. Finally the author is very grateful for the reviewers' very useful suggestions including the method for proving $w_2\not\leq w_1$ that is now used in subsections \ref{xadm16} and \ref{xadm17}.
\section{Notions}
\label{not}
Let us recall the notions of $\mu$-admissibility and $\mu$-permissibility introduced by Kottwitz and Rapoport in \cite{kotrapo}. 

Let $(X^{*},X_{*},R,R^{\vee})$ be a root datum with reduced root system $R$. Let $W$ be its (finite) Weyl group,  $W_a=\mathbb{Z}R^{\vee}\rtimes W$ its affine Weyl group and $\Wtilde=X_{*}\rtimes W$ its extended affine Weyl group.

Set $V:=X_{*}\otimes_{\mathbb{Z}}\mathbb{R}$. We choose a base $B=\lbrace\alpha_1,\dots, \alpha_l\rbrace$ of $R$ and denote by $\mathfrak{A}$ the base alcove
$$\mathfrak{A}=\lbrace v\in V\ \vert\  \langle \alpha_i,v\rangle> 0\text{ for }i=1,\dots,l\text{ and } \langle \tilde{\alpha},v\rangle< 1\rbrace$$
in $V$, where
$$\tilde{\alpha}=\sum_{i=1}^{l}n_i\alpha_i$$
is the highest root of $R$. We consider the stabilizer
$$\Omega:=\lbrace x\inWtilde \vert\  x(\mathfrak{A})=\mathfrak{A}\rbrace$$
of $\mathfrak{A}$ under the $\Wtilde$-operation on the alcoves in $V$. Every element $x\inWtilde $ can be written uniquely in the form $x=y\omega$ with $y\in W_a$ and $\omega\in \Omega$, since $W_a$ acts freely and transitively on the set of alcoves.

Now we can extend the Bruhat order on the Coxeter group $W_a$ to $\Wtilde$ as follows: For elements $x=y\omega$ and $x'=y'\omega'$ of the extended affine Weyl group $\Wtilde $ with $y,y'\in W_a$ and $\omega, \omega'\in \Omega$, we set $x\leq x'$, if $\omega=\omega'$ and $y\leq y'$ in the Bruhat order on $W_a$.

If we consider some $\nu\in X_{*}$ as an element of $\Wtilde$, we denote it by $t_\nu$. Thus, $t_\nu$ is the translation $V\to V$, $v\mapsto v+\nu$.

Let $\mu\in X_{*}$ be a fixed cocharacter.

An element $x\inWtilde$ of the extended affine Weyl group is called \emph{$\mu$-admissible}, if $x\leq t_{w(\mu)}$ for some $w\in W$. The set of all $\mu$-admissible $x\inWtilde$ is denoted by $\Adm(\mu)$.

Let $P_{\mu}$ be the convex hull in $V$ of the $W$-orbit $\Wmu=\lbrace w(\mu)\ \vert\  w\in W\rbrace$. An element $x\inWtilde $ of the extended affine Weyl group is called \emph{$\mu$-permissible}, if it satisfies both of the following conditions:
\begin{itemize}
\item[(i)] If $x=y_1\omega_1$ and $t_{\mu}=y_2\omega_2$ with $y_1,y_2\in W_a$ and $\omega_1,\omega_2\in \Omega$, we have $\omega_1=\omega_2$.
\item[(ii)] For every element $v\in \overline{\mathfrak{A}}$ in the closure of the base alcove $\mathfrak{A}$ we have $x(v)-v\in P_\mu$.
\end{itemize}
The first condition (i) is equivalent to $x\in W_at_{\mu}$. The set of all $\mu$-permissible $x\inWtilde$ is denoted by $\Perm(\mu)$.

Let $\lbrace\rho_1,\dots, \rho_l\rbrace$ be the dual basis in $\mathbb{Z}R^{\vee}\otimes_{\mathbb{Z}}\mathbb{R}$ of the basis $\lbrace\alpha_1,\dots, \alpha_l\rbrace$ of $\mathbb{Z}R\otimes_{\mathbb{Z}}\mathbb{R}$. Furthermore, we set
$$a_1:=\frac{\rho_1}{n_1}, \dots, a_{l}:=\frac{\rho_l}{n_l}\text{ and }a_{l+1}:=0.$$
Then $a_1$, \dots, $a_{l+1}$ are elements of the minimal facets of $\overline{\mathfrak{A}}$. Therefore, condition (ii) in the definition of $\mu$-permissibility is equivalent to  $x(a_i)-a_i\in P_\mu$ for $i=1,\dots,l+1$. In the following we will consider root data, where $a_1$, \dots, $a_{l+1}$ are in fact the vertices of $\overline{\mathfrak{A}}$ and every element of $\overline{\mathfrak{A}}$ is a convex combination of $a_1$, \dots, $a_{l+1}$. Here, the equivalence of condition (ii) and $x(a_i)-a_i\in P_\mu$ for $i=1,\dots,l+1$ is easy to see.

Recall that $\mu\in X_{*}$ is called minuscule, if $\langle \alpha, \mu\rangle\in \lbrace -1,0,1\rbrace$ for all $\alpha\in R$.
\section{The case of $E_6$}
\label{e6}
\subsection{Setting}
\label{settin6}
Let $R$ denote the root system $E_6$. Furthermore, let $R^{\vee}$ be the dual root system, $X^{*}:=Q(R)=\mathbb{Z}R$ the root lattice and $X_{*}:=P(R^{\vee})$ the coweight lattice. Then $(X^{*},X_{*},R,R^{\vee})$ is a root datum with (finite) Weyl group $W$, affine Weyl group  $W_a=Q(R^{\vee})\rtimes W$ and extended affine Weyl group $\Wtilde=X_{*}\rtimes W=P(R^{\vee})\rtimes W$. Observe that by definition the cocharacters $X_{*}$ agree with the coweights $P(R^{\vee})$.

For our calculations we consider the explicit construction of $R$ given in \cite[Plate V(I)]{bour}: Let
$$V^{*}:=\lbrace (x_1,\dots,x_8)\in \mathbb{R}^{8}\ \vert\ x_6=x_7=-x_8\rbrace.$$
Furthermore, let $e_1,\dots, e_8$ be the standard basis vectors of $\mathbb{R}^{8}$. Then the 72 vectors in $V^{*}$
$$\pm e_i\pm e_j$$
for $1 \leq i<j\leq 5$ and
$$\pm \frac{1}{2}\left(e_8-e_7-e_6+\sum_{i=1}^{5}(-1)^{\delta(i)}e_i\right)$$
with $\sum_{i=1}^{5} \delta(i)$ even form a root system of type $E_6$, which we identify with $R$.

Let $e_1',\dots, e_8'$ be the dual basis of $e_1,\dots, e_8$ in the dual space of $\mathbb{R}^{8}$. The dual space of $V^{*}$ can be identified with
$$V:=\lbrace x_1'e_1'+\dots+x_8'e_8'\in \mathbb{R}^{8}\ \vert\ x_6'=x_7'=-x_8'\rbrace$$
and the dual root system $R^{\vee}$ consists of
$$\pm e_i'\pm e_j'$$
for $1 \leq i<j\leq 5$ and
$$\pm \frac{1}{2}\left(e_8'-e_7'-e_6'+\sum_{i=1}^{5}(-1)^{\delta(i)}e_i'\right)$$
with $\sum_{i=1}^{5} \delta(i)$ even, see \cite[Plate V(V)]{bour}.

Then $V^{*}=Q(R)\otimes_{\mathbb{Z}}\mathbb{R}=X^{*}\otimes_{\mathbb{Z}}\mathbb{R}$ and $V=Q(R^{\vee})\otimes_{\mathbb{Z}}\mathbb{R}=P(R^{\vee})\otimes_{\mathbb{Z}}\mathbb{R}=X_{*}\otimes_{\mathbb{Z}}\mathbb{R}$, which agrees with the notation from section \ref{not}.

As in \cite[Plate V(II)]{bour} we consider the base of $R$ given by
\begin{eqnarray*}
\alpha_1&=&\frac{1}{2}(e_1+e_8)-\frac{1}{2}(e_2+e_3+e_4+e_5+e_6+e_7)\\
\alpha_2&=&e_1+e_2\\
\alpha_3&=&e_2-e_1\\
\alpha_4&=&e_3-e_2\\
\alpha_5&=&e_4-e_3\\
\alpha_6&=&e_5-e_4.
\end{eqnarray*}
The Dynkin diagram is (see \cite[Plate V(IV)]{bour}):

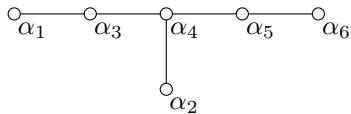
\begin{center}
\begin{tikzpicture}[line cap=round,line join=round,>=triangle 45,x=0.5cm,y=0.5cm]
\clip(-4,-1) rectangle (6,3);
\draw(-3,2) circle (0.08cm);
\draw(-1,2) circle (0.08cm);
\draw(1,2) circle (0.08cm);
\draw(3,2) circle (0.08cm);
\draw(5,2) circle (0.08cm);
\draw(1,0) circle (0.08cm);
\draw (-2.84,2)-- (-1.16,2);
\draw (0.84,2)-- (-0.84,2);
\draw (1.16,2)-- (2.84,2);
\draw (4.84,2)-- (3.16,2);
\draw (1,1.84)-- (1,0.16);
\draw (-3.15,2.02) node[anchor=north west] {$ \alpha_1 $};
\draw (-1.15,2.02) node[anchor=north west] {$ \alpha_3 $};
\draw (0.85,2.02) node[anchor=north west] {$ \alpha_4 $};
\draw (2.85,2.02) node[anchor=north west] {$ \alpha_5 $};
\draw (4.85,2.02) node[anchor=north west] {$ \alpha_6 $};
\draw (0.84,0.02) node[anchor=north west] {$ \alpha_2 $};
\end{tikzpicture}
\end{center}

The simple reflections belonging to $\alpha_1$, \dots, $\alpha_6$ are as usual denoted by $s_1$, \dots, $s_6$, respectively. Thus, $s_i: V\to V$ is given by
$$s_i(v)=v-\langle \alpha_i, v\rangle \alpha_i^{\vee}.$$

According to \cite[Plate V(IV)]{bour} the highest root of $R$ is $\tilde{\alpha}=\alpha_1+2\alpha_2+2\alpha_3+3\alpha_4+2\alpha_5+\alpha_6$. Thus, 
$$\mu:=\rho_1=\frac{2}{3}(e_8'-e_7'-e_6')$$
is a dominant minuscule coweight (see \cite[Plate V(VI)]{bour}).

We set $x=w_2t_\mu w_1^{-1}\inWtilde$ with
$$w_1=s_{2}s_{4}s_{5}s_{6}s_{3}s_{4}s_{5}s_{2}s_{4}s_{3}s_{1}$$
and
$$w_2=s_{4}s_{5}s_{6}s_{2}s_{4}s_{5}.$$

In the following subsections we will check that $x\in\Perm(\mu)$, but $x\not\in\Adm(\mu)$. Hence $\Perm(\mu)\neq\Adm(\mu)$ for $\mu=\rho_1$. As there is an automorphism of the root system $E_6$ interchanging $\alpha_1$ and $\alpha_6$,  we also get $\Perm(\rho_6)\neq\Adm(\rho_6)$. Note that $\rho_1$ and $\rho_6$ are the only dominant minuscule coweights for $E_6$.
\subsection{$x\in\Perm(\mu)$}
\label{xperm6}
In this subsection we prove that $x$ is $\mu$-permissible. The orbit $\Wmu$ consists of the following 27 elements:
$$\mu=\frac{2}{3}(e_8'-e_7'-e_6'),$$
$$\frac{1}{6}(e_8'-e_7'-e_6')-\frac{1}{2}\left(\sum_{i=1}^{5}(-1)^{\delta(i)}e_i'\right)$$
with $\sum_{i=1}^{5} \delta(i)$ even and
$$-\frac{1}{3}(e_8'-e_7'-e_6')\pm e_i'$$
with $1\leq i\leq 5$.

Now $P_\mu$ is the convex hull of these 27 points and we have to show $x(a_i)-a_i\in P_\mu$ for $i=1,\dots, 7$. Here, the $a_i$ are given by (see \cite[Plate V(IV)]{bour})
$$a_1=\frac{\rho_1}{1}=\mu=\begin{pmatrix}0\\ 0\\ 0\\ 0\\ 0\\ -\frac{2}{3}\\ -\frac{2}{3}\\ \frac{2}{3}\end{pmatrix},
\quad a_2=\frac{\rho_2}{2}=\begin{pmatrix}\frac{1}{4}\\ \frac{1}{4}\\ \frac{1}{4}\\ \frac{1}{4}\\ \frac{1}{4}\\ -\frac{1}{4}\\ -\frac{1}{4}\\ \frac{1}{4}\end{pmatrix},
\quad a_3=\frac{\rho_3}{2}=\begin{pmatrix}-\frac{1}{4}\\ \frac{1}{4}\\ \frac{1}{4}\\ \frac{1}{4}\\ \frac{1}{4}\\ -\frac{5}{12}\\ -\frac{5}{12}\\ \frac{5}{12}\end{pmatrix},
\quad a_4=\frac{\rho_4}{3}=\begin{pmatrix}0\\ 0\\ \frac{1}{3}\\ \frac{1}{3}\\ \frac{1}{3}\\ -\frac{1}{3}\\ -\frac{1}{3}\\ \frac{1}{3}\end{pmatrix},$$
$$a_5=\frac{\rho_5}{2}=\begin{pmatrix}0\\ 0\\ 0\\ \frac{1}{2}\\ \frac{1}{2}\\ -\frac{1}{3}\\ -\frac{1}{3}\\ \frac{1}{3}\end{pmatrix},
\quad a_6=\frac{\rho_6}{1}=\begin{pmatrix}0\\ 0\\ 0\\ 0\\ 1\\ -\frac{1}{3}\\ -\frac{1}{3}\\ \frac{1}{3}\end{pmatrix},
\quad a_7=0=\begin{pmatrix}0\\ 0\\ 0\\ 0\\ 0\\ 0\\ 0\\ 0\end{pmatrix}.$$

Since $\mu$ is fixed by $s_2$, $s_3$, $s_4$, $s_5$ and $s_6$, we have
$$x=w_2t_\mu w_1^{-1}=t_{w_2(\mu)}w_2w_1^{-1}=t_\mu w_2w_1^{-1}=t_\mu s_{4}s_{5}s_{6}s_{2}s_{4}s_{5}s_{1}s_{3}s_{4}s_{2}s_{5}s_{4}s_{3}s_{6}s_{5}s_{4}s_{2}.$$
Multiplying the matrices for the simple reflections $s_i$ we get that $w_2w_1^{-1}$ is represented by the matrix
$$M:=\frac{1}{4}\begin{pmatrix}
-1&3&-1&-1&-1&1&1&-1\\
-3&1&1&1&1&-1&-1&1\\
1&1&1&1&-3&-1&-1&1\\
1&1&1&-3&1&-1&-1&1\\
1&1&-3&1&1&-1&-1&1\\
1&1&1&1&1&3&-1&1\\
1&1&1&1&1&-1&3&1\\
-1&-1&-1&-1&-1&1&1&3
\end{pmatrix}.$$

Now $x(a_i)-a_i=M a_i+\mu-a_i$ for $i=1,\dots,7$ yields the following results
$$x(a_1)-a_1=\begin{pmatrix} -\frac{1}{2}\\ \frac{1}{2}\\ \frac{1}{2}\\ \frac{1}{2}\\ \frac{1}{2}\\ -\frac{1}{6}\\ -\frac{1}{6}\\ \frac{1}{6} \end{pmatrix}, \quad\quad
x(a_2)-a_2=\begin{pmatrix} -\frac{1}{2}\\ 0\\ 0\\ 0\\ 0\\ -\frac{1}{6}\\ -\frac{1}{6}\\ \frac{1}{6} \end{pmatrix}=\frac{1}{2}\begin{pmatrix}-\frac{1}{2}\\ \frac{1}{2}\\ \frac{1}{2}\\ \frac{1}{2}\\ \frac{1}{2}\\ -\frac{1}{6}\\ -\frac{1}{6}\\ \frac{1}{6} \end{pmatrix}+\frac{1}{2}\begin{pmatrix}-\frac{1}{2}\\ -\frac{1}{2}\\ -\frac{1}{2}\\ -\frac{1}{2}\\ -\frac{1}{2}\\ -\frac{1}{6}\\ -\frac{1}{6}\\ \frac{1}{6} \end{pmatrix},$$
$$x(a_3)-a_3=\begin{pmatrix} 0\\ \frac{1}{2}\\ 0\\ 0\\ 0\\ -\frac{1}{6}\\ -\frac{1}{6}\\ \frac{1}{6} \end{pmatrix}=\frac{1}{2}\begin{pmatrix}-\frac{1}{2}\\ \frac{1}{2}\\ \frac{1}{2}\\ \frac{1}{2}\\ \frac{1}{2}\\ -\frac{1}{6}\\ -\frac{1}{6}\\ \frac{1}{6} \end{pmatrix}+\frac{1}{2}\begin{pmatrix}\frac{1}{2}\\ \frac{1}{2}\\ -\frac{1}{2}\\ -\frac{1}{2}\\ -\frac{1}{2}\\ -\frac{1}{6}\\ -\frac{1}{6}\\ \frac{1}{6} \end{pmatrix},$$
$$x(a_4)-a_4=\begin{pmatrix} -\frac{1}{2}\\ \frac{1}{2}\\ -\frac{1}{6}\\ -\frac{1}{6}\\ -\frac{1}{6}\\ -\frac{1}{6}\\ -\frac{1}{6}\\ \frac{1}{6}\end{pmatrix}=\frac{1}{3}\begin{pmatrix}-\frac{1}{2}\\ \frac{1}{2}\\ -\frac{1}{2}\\ -\frac{1}{2}\\ \frac{1}{2}\\ -\frac{1}{6}\\ -\frac{1}{6}\\ \frac{1}{6} \end{pmatrix}+\frac{1}{3}\begin{pmatrix}-\frac{1}{2}\\ \frac{1}{2}\\ -\frac{1}{2}\\ \frac{1}{2}\\ -\frac{1}{2}\\ -\frac{1}{6}\\ -\frac{1}{6}\\ \frac{1}{6} \end{pmatrix}+\frac{1}{3}\begin{pmatrix} -\frac{1}{2}\\ \frac{1}{2}\\ \frac{1}{2}\\ -\frac{1}{2}\\ -\frac{1}{2}\\ -\frac{1}{6}\\ -\frac{1}{6}\\ \frac{1}{6} \end{pmatrix},$$
$$x(a_5)-a_5=\begin{pmatrix} -\frac{1}{2}\\ \frac{1}{2}\\ 0\\ -\frac{1}{2}\\ 0\\ -\frac{1}{6}\\ -\frac{1}{6}\\ \frac{1}{6} \end{pmatrix}=\frac{1}{2}\begin{pmatrix}-\frac{1}{2}\\ \frac{1}{2}\\ -\frac{1}{2}\\ -\frac{1}{2}\\ \frac{1}{2}\\ -\frac{1}{6}\\ -\frac{1}{6}\\ \frac{1}{6} \end{pmatrix}+\frac{1}{2}\begin{pmatrix}-\frac{1}{2}\\ \frac{1}{2}\\ \frac{1}{2}\\ -\frac{1}{2}\\ -\frac{1}{2}\\ -\frac{1}{6}\\ -\frac{1}{6}\\ \frac{1}{6} \end{pmatrix},$$
$$x(a_6)-a_6=\begin{pmatrix} -\frac{1}{2}\\ \frac{1}{2}\\ -\frac{1}{2}\\ \frac{1}{2}\\ -\frac{1}{2}\\ -\frac{1}{6}\\ -\frac{1}{6}\\ \frac{1}{6} \end{pmatrix},\quad\quad
x(a_7)-a_7=\begin{pmatrix} 0\\ 0\\ 0\\ 0\\ 0\\ -\frac{2}{3}\\ -\frac{2}{3}\\ \frac{2}{3} \end{pmatrix}.$$
Thus, every $x(a_i)-a_i$ for $i=1,\dots,7$ is a convex combination of points in $W\mu$ and therefore $x(a_i)-a_i\in P_\mu$. Hence, condition (ii) for $x\in \Perm (\mu)$ is fulfilled. On the other hand condition (i) is trivial:
$$x=w_2t_\mu w_1^{-1}=w_2w_1^{-1}t_{w_1(\mu)}\in W_at_{w_1(\mu)}=W_at_{\mu}.$$

All calculations in this subsection were both done by hand and with the computer software Octave.
\subsection{First test of $x\not\in\Adm(\mu)$}
\label{xadm16}
Let us denote by
$$I(\mu):=\lbrace s_i\ \vert \ i\in \lbrace 1,\dots, l\rbrace,\langle \alpha_i,\mu\rangle=0\rbrace$$
the set of simple reflections fixing $\mu$. Since $\mu=\rho_1$, we have $I(\mu)=\lbrace s_2,s_3,s_4,s_5,s_6\rbrace$ in our case. 

Now let $W_{I(\mu)}$ be the subgroup of $W$ generated by $I(\mu)$. Consider the set
$$W^{I(\mu)}=\lbrace w\in W\ \vert\  w\leq ws \text{ for all }s\in I(\mu)\rbrace$$
of minimal length representatives in $W/W_{I(\mu)}$.

According to He and Lam \cite[Theorem 2.2(1)]{lamhe} the map
$$W^{I(\mu)}\times W\to Wt_\mu W,\quad (z_1,z_2)\mapsto z_2t_\mu z_1^{-1}.$$
is  a bijection (see also \cite[Proposition 3.2(1)]{he}). Furthermore \cite[Theorem 2.2(2)]{lamhe} states that an element $z_2t_\mu z_1^{-1}\in \Wtilde$ with $z_1\in W^{I(\mu)}$ and $z_2\in W$ is $\mu$-admissible iff $z_2\leq z_1$ (see also \cite[Proposition 3.2(2)]{he}).

Thus, for proving $x=w_2t_\mu w_1^{-1}\not\in\Adm(\mu)$ it is sufficient to show $w_1\in W^{I(\mu)}$ and $w_2\not\leq w_1$. Recall that for $w\in W$ and a simple reflection $s_i$ (belonging to the root $\alpha_i$) we have $w\leq ws_i$ if $w(\alpha_i)$ is a positive root and $ws_i\leq w$ if $w(\alpha_i)$ is a negative root (see \cite[Theorem 5.4]{humph}).

For checking $w_2\not\leq w_1$ we will repeatedly use the following fact (which holds for Coxeter groups in general): \textit{For $w', w\in W$ with $w'\leq w$ and a simple reflection $s$ we have $w'\leq ws$ or $w's\leq ws$, so in particular $\min(w',w's)\leq ws$.}

\textit{Proof.} Let $ws=t_1\dots t_k$ be a reduced expression. If $w\leq ws$, then $w'\leq w\leq ws$ and we are done. Hence we may assume $ws\leq w$, then $w=t_1\dots t_ks$ is a reduced expression. By \cite[Theorem 5.10]{humph} $w'$ has a reduced expression being a subexpression of $t_1\dots t_ks$. If this subexpression does not contain the last $s$, then $w'\leq t_1\dots t_k=ws$. If the subexpression contains the last $s$, then $w's$ has a reduced expression being a subexpression of $t_1\dots t_k$ (namely the subexpression of $t_1\dots t_ks$ with the last $s$ dropped), hence $w's\leq t_1\dots t_k=ws$.\hfill\qedsymbol

Now let us assume for contradiction that $w_2\leq w_1$, i.e.
$$s_{4}s_{5}s_{6}s_{2}s_{4}s_{5}\leq s_{2}s_{4}s_{5}s_{6}s_{3}s_{4}s_{5}s_{2}s_{4}s_{3}s_{1}.$$
As $s_{4}s_{5}s_{6}s_{2}s_{4}s_{5}(\alpha_1)=\alpha_1$ is a positive root, we have
$$\min (s_{4}s_{5}s_{6}s_{2}s_{4}s_{5}, s_{4}s_{5}s_{6}s_{2}s_{4}s_{5}s_{1})=s_{4}s_{5}s_{6}s_{2}s_{4}s_{5}$$
and we can deduce
$$s_{4}s_{5}s_{6}s_{2}s_{4}s_{5}\leq s_{2}s_{4}s_{5}s_{6}s_{3}s_{4}s_{5}s_{2}s_{4}s_{3}.$$
Then by $s_{4}s_{5}s_{6}s_{2}s_{4}s_{5}(\alpha_3)=\alpha_2+\alpha_3+2\alpha_4+\alpha_5$ being a positive root we get
$$s_{4}s_{5}s_{6}s_{2}s_{4}s_{5}\leq s_{2}s_{4}s_{5}s_{6}s_{3}s_{4}s_{5}s_{2}s_{4}.$$
Repeating this process we get the following (including the initial steps explained above):

\begin{math}
\begin{array}{lcl}
&&s_{4}s_{5}s_{6}s_{2}s_{4}s_{5}\leq s_{2}s_{4}s_{5}s_{6}s_{3}s_{4}s_{5}s_{2}s_{4}s_{3}s_{1}\\
s_{4}s_{5}s_{6}s_{2}s_{4}s_{5}(\alpha_1)=\alpha_1>0&\quad \rightarrow \quad &s_{4}s_{5}s_{6}s_{2}s_{4}s_{5}\leq s_{2}s_{4}s_{5}s_{6}s_{3}s_{4}s_{5}s_{2}s_{4}s_{3}\\
s_{4}s_{5}s_{6}s_{2}s_{4}s_{5}(\alpha_3)=\alpha_2+\alpha_3+2\alpha_4+\alpha_5>0&\quad \rightarrow \quad &s_{4}s_{5}s_{6}s_{2}s_{4}s_{5}\leq s_{2}s_{4}s_{5}s_{6}s_{3}s_{4}s_{5}s_{2}s_{4}\\
s_{4}s_{5}s_{6}s_{2}s_{4}s_{5}(\alpha_4)=\alpha_6>0&\quad \rightarrow \quad &s_{4}s_{5}s_{6}s_{2}s_{4}s_{5}\leq s_{2}s_{4}s_{5}s_{6}s_{3}s_{4}s_{5}s_{2}\\
s_{4}s_{5}s_{6}s_{2}s_{4}s_{5}(\alpha_2)=\alpha_5>0&\quad \rightarrow \quad &s_{4}s_{5}s_{6}s_{2}s_{4}s_{5}\leq s_{2}s_{4}s_{5}s_{6}s_{3}s_{4}s_{5}\\
s_{4}s_{5}s_{6}s_{2}s_{4}s_{5}(\alpha_5)=-(\alpha_2+\alpha_4+\alpha_5+\alpha_6)<0&\quad \rightarrow \quad &s_{4}s_{5}s_{6}s_{2}s_{4}\leq s_{2}s_{4}s_{5}s_{6}s_{3}s_{4}\\
s_{4}s_{5}s_{6}s_{2}s_{4}(\alpha_4)=-(\alpha_2+\alpha_4+\alpha_5)<0&\quad \rightarrow \quad &s_{4}s_{5}s_{6}s_{2}\leq s_{2}s_{4}s_{5}s_{6}s_{3}\\
s_{4}s_{5}s_{6}s_{2}(\alpha_3)=\alpha_3+\alpha_4>0&\quad \rightarrow \quad &s_{4}s_{5}s_{6}s_{2}\leq s_{2}s_{4}s_{5}s_{6}\\
s_{4}s_{5}s_{6}s_{2}(\alpha_6)=-(\alpha_4+\alpha_5+\alpha_6)<0&\quad \rightarrow \quad &s_{4}s_{5}s_{2}=s_{4}s_{5}s_{6}s_{2}s_{6}\leq s_{2}s_{4}s_{5}\\
s_{4}s_{5}s_{2}(\alpha_5)=-(\alpha_4+\alpha_5)<0&\quad \rightarrow \quad &s_{4}s_{2}=s_{4}s_{5}s_{2}s_{5}\leq s_{2}s_{4}\\
s_{4}s_{2}(\alpha_4)=\alpha_2>0&\quad \rightarrow \quad &s_{4}s_{2}\leq s_{2}
\end{array}
\end{math}

So we get $s_{4}s_{2}\leq s_{2}$, which is a contradiction. Hence we must have $w_2\not\leq w_1$.

For showing
$$w_1\in W^{I(\mu)}=\lbrace w\in W\ \vert\  w\leq ws_i \text{ for }i=2,3,4,5,6\rbrace$$
we just have to check that $w_1(\alpha_i)$ is a positive root for $i=2,3,4,5,6$. Indeed
\begin{eqnarray*}
w_1(\alpha_2)&=&\alpha_6>0\\
w_1(\alpha_3)&=&\alpha_1+\alpha_2+\alpha_3+\alpha_4>0\\
w_1(\alpha_4)&=&\alpha_5>0\\
w_1(\alpha_5)&=&\alpha_4>0\\
w_1(\alpha_6)&=&\alpha_3>0.
\end{eqnarray*}
All calculations above have been carried out by hand and were afterwards double-checked using Octave. Furthermore the statements $w_2\not\leq w_1$ and $w_1\in W^{I(\mu)}$ were also verified with the computer software Sage.
\subsection{Second test of $x\not\in\Adm(\mu)$}
\label{xadm26}
Using the computer software Sage the statement $x\not\in\Adm(\mu)$ is checked again, this time relying on a result of Haines instead of the result by He and Lam used in the last subsection. In \cite[proof of Proposition 4.6]{haines} Haines showed $w\leq t_{w(0)}$ for every $\mu$-admissible element $w\inWtilde$. Thus, for proving $x\not\in\Adm(\mu)$ it is enough to check $x\not\leq t_{\mu}$ (as $x(0)=w_2t_\mu w_1^{-1}(0)=w_2(\mu)=\mu$).

Sage provides a function \verb+reduced_word_of_alcove_morphism+, that calculates for any $x'\inWtilde$ the corresponding element $y'\in W_a$ such that $x'=y'\omega$ for some $\omega\in \Omega$.

The following calculation works with weights instead of coweights, but this does not make a difference since the root system $E_6$ is self-dual.

\begin{verbatim}
sage: R = RootSystem(["E",6,1]).weight_lattice()
sage: Lambda = R.fundamental_weights()
sage: W = WeylGroup(R)
sage: s = W.simple_reflections()
sage: R.reduced_word_of_alcove_morphism((Lambda[1]-Lambda[0]).translation)
[0, 2, 4, 3, 5, 4, 2, 0, 6, 5, 4, 2, 3, 4, 5, 6]
sage: y1=s[0]*s[2]*s[4]*s[3]*s[5]*s[4]*s[2]*s[0]*s[6]*s[5]*s[4]*s[2]*s[3]*s[4]*s[5]*s[6]
sage: w1=s[2]*s[4]*s[5]*s[6]*s[3]*s[4]*s[5]*s[2]*s[4]*s[3]*s[1]
sage: w1inv=s[1]*s[3]*s[4]*s[2]*s[5]*s[4]*s[3]*s[6]*s[5]*s[4]*s[2]
sage: w2=s[4]*s[5]*s[6]*s[2]*s[4]*s[5]
sage: w1mu=w1.action(Lambda[1]-Lambda[0])
sage: R.reduced_word_of_alcove_morphism(w1mu.translation)
[2, 4, 3, 5, 4, 2, 0, 6, 5, 4, 2, 3, 1, 4, 5, 6]
sage: y2=w2*w1inv*s[2]*s[4]*s[3]*s[5]*s[4]*s[2]*s[0]*s[6]*s[5]*s[4]*s[2]*s[3]*s[1]*s[4]*s[5]*s[6]
sage: y2.bruhat_le(y1)
False
\end{verbatim}

Here \verb+R+ is defined to be the affine root system $E_6$, \verb+Lambda+ denotes the fundamental weights, \verb+W+ is the affine Weyl group $W_a$ and \verb+s+ denotes the simple reflections. In the 5th line the function \verb+reduced_word_of_alcove_morphism+ computes a reduced word for $y_1\in W_a$ such that $t_{\mu}=y_1\omega_1$ for some $\omega_1\in \Omega$ (recall that $\mu=\rho_1$) and \verb+y1+ is this element $y_1\in W_a$. Then \verb+w1+, \verb+w1inv+ and \verb+w2+ are defined as $w_1$, $w_1^{-1}$ and $w_2$ in subsection \ref{settin6}, respectively. Furthermore, \verb+w1mu+ is defined to be $w_1(\mu)$. In the 11th line the function \verb+reduced_word_of_alcove_morphism+ computes a reduced word for $y'\in W_a$ such that $t_{w_1(\mu)}=y'\omega_2$ for some $\omega_2\in \Omega$. The next line defines \verb+y2+ to be $y_2:=w_2w_1^{-1}y'$. In fact, we have
$$x=w_2t_\mu w_1^{-1}=w_2w_1^{-1}t_{w_1(\mu)}=w_2w_1^{-1}y'\omega_2=y_2\omega_2$$
with $y_2\in W_a$ and $\omega_2\in \Omega$. For checking $x\not\leq t_{\mu}$ it is therefore enough to ask whether $y_2\leq y_1$. This is done in the last line. Indeed, Sage returns \verb+False+.

\section{The case of $E_7$}
\label{e7}
This section is basically a repetition of the last section, this time working with $E_7$ instead of $E_6$.
\subsection{Setting}
\label{settin7}
Let $R$ now denote the root system $E_7$. As before, let $R^{\vee}$ be the dual root system, $X^{*}:=Q(R)=\mathbb{Z}R$ the root lattice and $X_{*}:=P(R^{\vee})$ the coweight lattice. Then $(X^{*},X_{*},R,R^{\vee})$ is a root datum with (finite) Weyl group $W$, affine Weyl group  $W_a=Q(R^{\vee})\rtimes W$ and extended affine Weyl group $\Wtilde=X_{*}\rtimes W=P(R^{\vee})\rtimes W$. Observe that by definition the cocharacters $X_{*}$ agree with the coweights $P(R^{\vee})$.

For our calculations we consider the explicit construction of $R$ given in \cite[Plate VI(I)]{bour}: Let
$$V^{*}:=\lbrace (x_1,\dots,x_8)\in \mathbb{R}^{8}\ \vert\ x_7=-x_8\rbrace.$$
Furthermore, let $e_1,\dots, e_8$ be the standard basis vectors of $\mathbb{R}^{8}$. Then the 126 vectors in $V^{*}$
$$\pm e_i\pm e_j$$
for $1 \leq i<j\leq 6$,
$$\pm (e_7-e_8)$$
and
$$\pm \frac{1}{2}\left(e_7-e_8+\sum_{i=1}^{6}(-1)^{\delta(i)}e_i\right)$$
with $\sum_{i=1}^{6} \delta(i)$ odd form a root system of type $E_7$, which we identify with $R$.

Let $e_1',\dots, e_8'$ be the dual basis of $e_1,\dots, e_8$ in the dual space of $\mathbb{R}^{8}$. The dual space of $V^{*}$ can be identified with
$$V:=\lbrace x_1'e_1'+\dots+x_8'e_8'\in \mathbb{R}^{8}\ \vert\ x_7'=-x_8'\rbrace$$
and the dual root system $R^{\vee}$ consists of
$$\pm e_i'\pm e_j'$$
for $1 \leq i<j\leq 6$,
$$\pm (e_7'-e_8')$$
and
$$\pm \frac{1}{2}\left(e_7'-e_8'+\sum_{i=1}^{6}(-1)^{\delta(i)}e_i'\right)$$
with $\sum_{i=1}^{6} \delta(i)$ odd, see \cite[Plate VI(V)]{bour}.

Then $V^{*}=Q(R)\otimes_{\mathbb{Z}}\mathbb{R}=X^{*}\otimes_{\mathbb{Z}}\mathbb{R}$ and $V=Q(R^{\vee})\otimes_{\mathbb{Z}}\mathbb{R}=P(R^{\vee})\otimes_{\mathbb{Z}}\mathbb{R}=X_{*}\otimes_{\mathbb{Z}}\mathbb{R}$, which agrees with the notation from section \ref{not}.

As in \cite[Plate VI(II)]{bour} we consider the base of $R$ given by
\begin{eqnarray*}
\alpha_1&=&\frac{1}{2}(e_1+e_8)-\frac{1}{2}(e_2+e_3+e_4+e_5+e_6+e_7)\\
\alpha_2&=&e_1+e_2\\
\alpha_3&=&e_2-e_1\\
\alpha_4&=&e_3-e_2\\
\alpha_5&=&e_4-e_3\\
\alpha_6&=&e_5-e_4\\
\alpha_7&=&e_6-e_5
\end{eqnarray*}
The Dynkin diagram is (see \cite[Plate VI(IV)]{bour}):

\begin{center}
\begin{tikzpicture}[line cap=round,line join=round,>=triangle 45,x=0.5cm,y=0.5cm]
\clip(-4,-1) rectangle (8,3);
\draw(-3,2) circle (0.08cm);
\draw(-1,2) circle (0.08cm);
\draw(1,2) circle (0.08cm);
\draw(3,2) circle (0.08cm);
\draw(5,2) circle (0.08cm);
\draw(7,2) circle (0.08cm);
\draw(1,0) circle (0.08cm);
\draw (-2.84,2)-- (-1.16,2);
\draw (0.84,2)-- (-0.84,2);
\draw (1.16,2)-- (2.84,2);
\draw (4.84,2)-- (3.16,2);
\draw (6.84,2)-- (5.16,2);
\draw (1,1.84)-- (1,0.16);
\draw (-3.15,2.02) node[anchor=north west] {$ \alpha_1 $};
\draw (-1.15,2.02) node[anchor=north west] {$ \alpha_3 $};
\draw (0.85,2.02) node[anchor=north west] {$ \alpha_4 $};
\draw (2.85,2.02) node[anchor=north west] {$ \alpha_5 $};
\draw (4.85,2.02) node[anchor=north west] {$ \alpha_6 $};
\draw (6.85,2.02) node[anchor=north west] {$ \alpha_7 $};
\draw (0.84,0.02) node[anchor=north west] {$ \alpha_2 $};
\end{tikzpicture}
\end{center}

The simple reflections belonging to $\alpha_1$, \dots, $\alpha_7$ are as usual denoted by $s_1$, \dots, $s_7$, respectively. Thus, $s_i: V\to V$ is given by
$$s_i(v)=v-\langle \alpha_i, v\rangle \alpha_i^{\vee}.$$

According to \cite[Plate VI(IV)]{bour},  the highest root of $R$ is $\tilde{\alpha}=2\alpha_1+2\alpha_2+3\alpha_3+4\alpha_4+3\alpha_5+2\alpha_6+\alpha_7$. Thus, 
$$\mu:=\rho_7=\frac{1}{2}(e_8'-e_7')+e_6'$$
is the only dominant minuscule coweight (see \cite[Plate VI(VI)]{bour}).

We set $x=w_2t_\mu w_1^{-1}\inWtilde$ with
$$w_1=s_{2}s_{4}s_{5}s_{3}s_{4}s_{1}s_{3}s_{2}s_{4}s_{5}s_{6}s_{7}$$
and
$$w_2=s_{4}s_{3}s_{2}s_{4}s_{1}s_{3}.$$

In the following subsections we will check that $x\in\Perm(\mu)$, but $x\not\in\Adm(\mu)$. This implies that $\Perm(\mu)\neq\Adm(\mu)$ for $\mu=\rho_7$.
\subsection{$x\in\Perm(\mu)$}
\label{xperm7}
In this subsection we prove that $x$ is $\mu$-permissible. The orbit $\Wmu$ consists of the following 56 elements:
$$\pm\frac{1}{2}(e_8'-e_7')\pm e_i'$$
with $1\leq i\leq 6$ and
$$\frac{1}{2}\left(\sum_{i=1}^{6}(-1)^{\delta(i)}e_i'\right)$$
with $\sum_{i=1}^{6} \delta(i)$ even. 

Now $P_\mu$ is the convex hull of these 56 points and we have to show $x(a_i)-a_i\in P_\mu$ for $i=1,\dots, 8$. Here, the $a_i$ are given by (see \cite[Plate V(IV)]{bour})
$$a_1=\frac{\rho_1}{2}=\begin{pmatrix}0\\ 0\\ 0\\ 0\\ 0\\ 0\\ -\frac{1}{2}\\ \frac{1}{2}\end{pmatrix},
\quad a_2=\frac{\rho_2}{2}=\begin{pmatrix}\frac{1}{4}\\ \frac{1}{4}\\ \frac{1}{4}\\ \frac{1}{4}\\ \frac{1}{4}\\ \frac{1}{4}\\ -\frac{1}{2}\\ \frac{1}{2}\end{pmatrix},
\quad a_3=\frac{\rho_3}{3}=\begin{pmatrix}-\frac{1}{6}\\ \frac{1}{6}\\ \frac{1}{6}\\ \frac{1}{6}\\ \frac{1}{6}\\ \frac{1}{6}\\ -\frac{1}{2}\\ \frac{1}{2}\end{pmatrix},
\quad a_4=\frac{\rho_4}{4}=\begin{pmatrix}0\\ 0\\ \frac{1}{4}\\ \frac{1}{4}\\ \frac{1}{4}\\ \frac{1}{4}\\ -\frac{1}{2}\\ \frac{1}{2}\end{pmatrix},$$
$$a_5=\frac{\rho_5}{3}=\begin{pmatrix}0\\ 0\\ 0\\ \frac{1}{3}\\ \frac{1}{3}\\ \frac{1}{3}\\ -\frac{1}{2}\\ \frac{1}{2}\end{pmatrix},
\quad a_6=\frac{\rho_6}{2}=\begin{pmatrix}0\\ 0\\ 0\\ 0\\ \frac{1}{2}\\ \frac{1}{2}\\ -\frac{1}{2}\\ \frac{1}{2}\end{pmatrix},
\quad a_7=\frac{\rho_7}{1}=\mu=\begin{pmatrix}0\\ 0\\ 0\\ 0\\ 0\\ 1\\ -\frac{1}{2}\\ \frac{1}{2}\end{pmatrix},
\quad a_8=0=\begin{pmatrix}0\\ 0\\ 0\\ 0\\ 0\\ 0\\ 0\\ 0\end{pmatrix}.$$

Since $\mu$ is fixed by $s_1$, $s_2$, $s_3$, $s_4$ (and $s_5$, $s_6$), we have
$$x=w_2t_\mu w_1^{-1}=t_{w_2(\mu)}w_2w_1^{-1}=t_\mu w_2w_1^{-1}=t_\mu
s_{4}s_{3}s_{2}s_{4}s_{1}s_{3}s_{7}s_{6}s_{5}s_{4}s_{2}s_{3}s_{1}s_{4}s_{3}s_{5}s_{4}s_{2}.$$
Multiplying the matrices for the simple reflections $s_i$ we get that $w_2w_1^{-1}$ is represented by the matrix
$$M:=\frac{1}{4}\begin{pmatrix}
-1&1&-1&3&1&1&1&-1\\
1&-1&-3&1&-1&-1&-1&1\\
3&1&-1&-1&1&1&1&-1\\
1&-1&1&1&3&-1&-1&1\\
1&-1&1&1&-1&3&-1&1\\
-1&-3&-1&-1&1&1&1&-1\\
1&-1&1&1&-1&-1&3&1\\
-1&1&-1&-1&1&1&1&3
\end{pmatrix}.$$

Now $x(a_i)-a_i=M a_i+\mu-a_i$ for $i=1,\dots,7$ yields the following results
$$x(a_1)-a_1=\begin{pmatrix} -\frac{1}{4}\\ \frac{1}{4}\\ -\frac{1}{4}\\ \frac{1}{4}\\ \frac{1}{4}\\ \frac{3}{4}\\ -\frac{1}{4}\\ \frac{1}{4} \end{pmatrix}=
\frac{1}{2}\begin{pmatrix} -\frac{1}{2}\\ \frac{1}{2}\\ -\frac{1}{2}\\ \frac{1}{2}\\ \frac{1}{2}\\ \frac{1}{2}\\ 0\\ 0 \end{pmatrix}+\frac{1}{2}\begin{pmatrix} 0\\ 0\\ 0\\ 0\\ 0\\ 1\\ -\frac{1}{2}\\ \frac{1}{2} \end{pmatrix}, \quad\quad
x(a_2)-a_2=\begin{pmatrix} -\frac{1}{4}\\ -\frac{1}{4}\\ -\frac{1}{4}\\ \frac{1}{4}\\ \frac{1}{4}\\ \frac{1}{4}\\ -\frac{1}{4}\\ \frac{1}{4} \end{pmatrix}=\frac{1}{2}\begin{pmatrix} -\frac{1}{2}\\ -\frac{1}{2}\\ -\frac{1}{2}\\ \frac{1}{2}\\ \frac{1}{2}\\ -\frac{1}{2}\\ 0\\ 0 \end{pmatrix}+\frac{1}{2}\begin{pmatrix} 0\\ 0\\ 0\\ 0\\ 0\\ 1\\ -\frac{1}{2}\\ \frac{1}{2} \end{pmatrix},$$
$$x(a_3)-a_3=\begin{pmatrix} \frac{1}{6}\\ -\frac{1}{6}\\ -\frac{1}{2}\\ \frac{1}{6}\\ \frac{1}{6}\\ \frac{1}{2}\\ -\frac{1}{3}\\ \frac{1}{3} \end{pmatrix}
=\frac{1}{3}\begin{pmatrix} \frac{1}{2}\\ -\frac{1}{2}\\ -\frac{1}{2}\\ \frac{1}{2}\\ \frac{1}{2}\\ \frac{1}{2}\\ 0\\ 0 \end{pmatrix}
+\frac{1}{3}\begin{pmatrix} 0\\ 0\\ -1\\ 0\\ 0\\ 0\\ -\frac{1}{2}\\ \frac{1}{2} \end{pmatrix}
+\frac{1}{3}\begin{pmatrix} 0\\ 0\\ 0\\ 0\\ 0\\ 1\\ -\frac{1}{2}\\ \frac{1}{2} \end{pmatrix},$$
$$x(a_4)-a_4=\begin{pmatrix} 0\\ 0\\ -\frac{1}{2}\\ \frac{1}{4}\\ \frac{1}{4}\\ \frac{1}{2}\\ -\frac{1}{4}\\ \frac{1}{4}\end{pmatrix}
=\frac{1}{4}\begin{pmatrix} -\frac{1}{2}\\ \frac{1}{2}\\ -\frac{1}{2}\\ \frac{1}{2}\\ \frac{1}{2}\\ \frac{1}{2}\\ 0\\ 0 \end{pmatrix}
+\frac{1}{4}\begin{pmatrix} \frac{1}{2}\\ -\frac{1}{2}\\ -\frac{1}{2}\\ \frac{1}{2}\\ \frac{1}{2}\\ \frac{1}{2}\\ 0\\ 0 \end{pmatrix}
+\frac{1}{4}\begin{pmatrix} 0\\ 0\\ -1\\ 0\\ 0\\ 0\\ -\frac{1}{2}\\ \frac{1}{2} \end{pmatrix}
+\frac{1}{4}\begin{pmatrix} 0\\ 0\\ 0\\ 0\\ 0\\ 1\\ -\frac{1}{2}\\ \frac{1}{2} \end{pmatrix},$$
$$x(a_5)-a_5=\begin{pmatrix} \frac{1}{6}\\ \frac{1}{6}\\ -\frac{1}{6}\\ \frac{1}{6}\\ \frac{1}{6}\\ \frac{1}{2}\\ -\frac{1}{3}\\ \frac{1}{3} \end{pmatrix}
=\frac{1}{3}\begin{pmatrix} \frac{1}{2}\\ \frac{1}{2}\\ -\frac{1}{2}\\ \frac{1}{2}\\ -\frac{1}{2}\\ \frac{1}{2}\\ 0\\ 0 \end{pmatrix}
+\frac{1}{3}\begin{pmatrix} 0\\ 0\\ 0\\ 0\\ 1\\ 0\\ -\frac{1}{2}\\ \frac{1}{2} \end{pmatrix}
+\frac{1}{3}\begin{pmatrix} 0\\ 0\\ 0\\ 0\\ 0\\ 1\\ -\frac{1}{2}\\ \frac{1}{2} \end{pmatrix},$$
$$x(a_6)-a_6=\begin{pmatrix} 0\\ 0\\ 0\\ \frac{1}{2}\\ 0\\ \frac{1}{2}\\ -\frac{1}{2}\\ \frac{1}{2} \end{pmatrix}
=\frac{1}{2}\begin{pmatrix} 0\\ 0\\ 0\\ 1\\ 0\\ 0\\ -\frac{1}{2}\\ \frac{1}{2} \end{pmatrix}
+\frac{1}{2}\begin{pmatrix} 0\\ 0\\ 0\\ 0\\ 0\\ 1\\ -\frac{1}{2}\\ \frac{1}{2} \end{pmatrix},\quad\quad
x(a_7)-a_7=\begin{pmatrix} 0\\ 0\\ 0\\ 0\\ 1\\ 0\\ -\frac{1}{2}\\ \frac{1}{2} \end{pmatrix},\quad\quad
x(a_8)-a_8=\begin{pmatrix} 0\\ 0\\ 0\\ 0\\ 0\\ 1\\ -\frac{1}{2}\\ \frac{1}{2} \end{pmatrix}.$$
Thus, every $x(a_i)-a_i$ for $i=1,\dots,8$ is a convex combination of points in $W\mu$ and therefore $x(a_i)-a_i\in P_\mu$. Hence, condition (ii) for $x\in \Perm (\mu)$ is fulfilled. On the other hand,  condition (i) is trivial:
$$x=w_2t_\mu w_1^{-1}=w_2w_1^{-1}t_{w_1(\mu)}\in W_at_{w_1(\mu)}=W_at_{\mu}.$$

All calculations in this subsection were both done by hand and with the computer software Octave.
\subsection{$x\not\in\Adm(\mu)$}
\label{xadm17}
As mentioned before, all the calculations are completely analogous to the case of $E_6$. For more detailed explanations and comments the reader is referred to section \ref{e6} about $E_6$.
Note that, as $\mu=\rho_7$, we now have $I(\mu)=\lbrace s_1,s_2,s_3,s_4,s_5,s_6\rbrace$.

As the case of $E_6$ we assume for contradiction that $w_2\leq w_1$ and get the following:

\begin{math}
\begin{array}{lcl}
&&s_{4}s_{3}s_{2}s_{4}s_{1}s_{3}\leq s_{2}s_{4}s_{5}s_{3}s_{4}s_{1}s_{3}s_{2}s_{4}s_{5}s_{6}s_{7}\\
s_{4}s_{3}s_{2}s_{4}s_{1}s_{3}(\alpha_7)=\alpha_7>0&\quad \rightarrow \quad &s_{4}s_{3}s_{2}s_{4}s_{1}s_{3}\leq s_{2}s_{4}s_{5}s_{3}s_{4}s_{1}s_{3}s_{2}s_{4}s_{5}s_{6}\\
s_{4}s_{3}s_{2}s_{4}s_{1}s_{3}(\alpha_6)=\alpha_6>0&\quad \rightarrow \quad &s_{4}s_{3}s_{2}s_{4}s_{1}s_{3}\leq s_{2}s_{4}s_{5}s_{3}s_{4}s_{1}s_{3}s_{2}s_{4}s_{5}\\
s_{4}s_{3}s_{2}s_{4}s_{1}s_{3}(\alpha_5)=\alpha_2+\alpha_3+2\alpha_4+\alpha_5>0&\quad \rightarrow \quad &s_{4}s_{3}s_{2}s_{4}s_{1}s_{3}\leq s_{2}s_{4}s_{5}s_{3}s_{4}s_{1}s_{3}s_{2}s_{4}\\
s_{4}s_{3}s_{2}s_{4}s_{1}s_{3}(\alpha_4)=\alpha_1>0&\quad \rightarrow \quad &s_{4}s_{3}s_{2}s_{4}s_{1}s_{3}\leq s_{2}s_{4}s_{5}s_{3}s_{4}s_{1}s_{3}s_{2}\\
s_{4}s_{3}s_{2}s_{4}s_{1}s_{3}(\alpha_2)=\alpha_3>0&\quad \rightarrow \quad &s_{4}s_{3}s_{2}s_{4}s_{1}s_{3}\leq s_{2}s_{4}s_{5}s_{3}s_{4}s_{1}s_{3}\\
s_{4}s_{3}s_{2}s_{4}s_{1}s_{3}(\alpha_3)=-(\alpha_1+\alpha_2+\alpha_3+\alpha_4)<0&\quad \rightarrow \quad &s_{4}s_{3}s_{2}s_{4}s_{1}\leq s_{2}s_{4}s_{5}s_{3}s_{4}s_{1}\\
s_{4}s_{3}s_{2}s_{4}s_{1}(\alpha_1)=-(\alpha_1+\alpha_3+\alpha_4)<0&\quad \rightarrow \quad &s_{4}s_{3}s_{2}s_{4}\leq s_{2}s_{4}s_{5}s_{3}s_{4}\\
s_{4}s_{3}s_{2}s_{4}(\alpha_4)=-(\alpha_2+\alpha_3+\alpha_4)<0&\quad \rightarrow \quad &s_{4}s_{3}s_{2}\leq s_{2}s_{4}s_{5}s_{3}\\
s_{4}s_{3}s_{2}(\alpha_3)=-(\alpha_3+\alpha_4)<0&\quad \rightarrow \quad &s_{4}s_{2}=s_{4}s_{3}s_{2}s_{3}\leq s_{2}s_{4}s_{5}\\
s_{4}s_{2}(\alpha_5)=\alpha_4+\alpha_5>0&\quad \rightarrow \quad &s_{4}s_{2}\leq s_{2}s_{4}\\
s_{4}s_{2}(\alpha_4)=\alpha_2>0&\quad \rightarrow \quad &s_{4}s_{2}\leq s_{2}
\end{array}
\end{math}

Since $s_{4}s_{2}\leq s_{2}$ is a contradiction, we must have $w_2\not\leq w_1$.

For proving
$$w_1\in W^{I(\mu)}=\lbrace w\in W\ \vert\  w\leq ws_i \text{ for }i=1,2,3,4,5,6\rbrace$$
it suffices to see that the following are positive roots:
\begin{eqnarray*}
w_1(\alpha_1)&=&\alpha_5>0\\
w_1(\alpha_2)&=&\alpha_1>0\\
w_1(\alpha_3)&=&\alpha_4>0\\
w_1(\alpha_4)&=&\alpha_3>0\\
w_1(\alpha_5)&=&\alpha_2+\alpha_4+\alpha_5+\alpha_6>0\\
w_1(\alpha_6)&=&\alpha_7>0.
\end{eqnarray*}
Again, all calculations were carried out by hand and afterwards double-checked using Octave. Also, $w_2\not\leq w_1$ and $w_1\in W^{I(\mu)}$ were again verified with Sage.

The statements $w_2\not\leq w_1$ for $E_6$ in section \ref{e6} and for $E_7$ in this section can also be deduced from each other: If we apply the nontrivial automorphism of the root system $E_6$ and then embed the Weyl group of $E_6$ into the Weyl group of $E_7$ (which preserves the Bruhat order), the elements  $w_1^{E_6}$ and $w_2^{E_6}$ from section \ref{e6} turn into $w_1^{E_7}s_7$ and $w_2^{E_7}$ from this section. Hence we have $w_2^{E_6}\leq w_1^{E_6}$ if and only if $w_2^{E_7}\leq w_1^{E_7}s_{7}$. As $w_2^{E_7}$ does not contain $s_7$ we have $\min(w_2^{E_7},w_2^{E_7}s_7)=w_2^{E_7}$ and therefore  $w_2^{E_7}\leq w_1^{E_7}s_{7}$ holds if and only if $w_2^{E_7}\leq w_1^{E_7}$ (by the fact discussed in subsection \ref{xadm16}). Hence we get $w_2^{E_6}\leq w_1^{E_6}$ if and only if $w_2^{E_7}\leq w_1^{E_7}$.

Finally, computations with Sage analogous to those presented in subsection \ref{xadm26} were carried out in this case as well and verified $x\not\in\Adm(\mu)$ again.

\textsc{Stanford University Department of Mathematics, 450 Serra Mall, Building 380, Stanford, CA 94305-2125, USA}

\textit{E-mail address:} \texttt{lsauerma@stanford.edu}

\textcopyright 2016. This manuscript version is made available under the CC-BY-NC-ND 4.0 license\\
http://creativecommons.org/licenses/by-nc-nd/4.0/
\end{document}